\documentclass{article}
\usepackage{amssymb, amsmath, euscript}

\newcommand{\Mr}{\mathbb R}
\newcommand{\ME}{\mathbf E}
\newcommand{\MP}{\mathbf P}

\newcommand{\CF}{\mathcal F}
\newcommand{\CG}{\mathcal G}
\newcommand{\CE}{\mathcal E}
\newcommand{\CX}{\mathcal X}
\newcommand{\CM}{\mathcal M}
\newcommand{\CS}{\mathcal S}
\newcommand{\SF}{\overline{F}}

\newcommand{\pa}{\partial}

\newcommand{\eps}{\varepsilon}

\newcommand{\Oa}{\Omega}

\newcommand{\qi}{\varphi}

\newcommand{\bs}{\bfseries}
\newcommand{\al}{\alpha}
\newcommand{\be}{\beta}
\newcommand{\ga}{\gamma}
\newcommand{\Ga}{\Gamma}

\newcommand{\bd}{\bar{\delta}}
\newcommand{\la}{\lambda}

\newcommand{\IL}{\int\limits}

\newcommand{\tF}{\tilde{F}}

\DeclareMathOperator{\ess}{ess}

\DeclareMathOperator{\vraimax}{vraimax}
\DeclareMathOperator{\mC}{C}
\DeclareMathOperator{\mD}{D}

\begin{document}
\large
\begin{center}

{\bf\LARGE Tail Asymptotic of Sum and Product of Random 

Variables with Applications in the Theory of 

Extremes of Conditionally Gaussian Processes}

\bigskip

{\it Andrey Sarantsev. June 2010}

{\it University of Washington}

{\it Department of Mathematics, PhD Student}

{\it E-mail: ansa1989@uw.edu}

\end{center}

{\normalsize {\bf Abstract.} We consider two independent random variables with the given tail asymptotic (e.g. power or exponential). We find tail asymptotics for their sum and product. This is done by some cumbersome but purely technical computations and requires the use of the Laplace method for asymptotic of integrals. We also recall the results for asymptotic of $\MP\{\sup_{t \ge 0}(X(t) - ct^{\be}) > u\}$ as $u \to \infty$, where $X = (X(t), t \ge 0)$ is a self-similar locally stationary centered Gaussian process; and we find the asymptotic for the same probability after replacing the constant $c$ by a random variable $\eta$, independent of $X$. We also find the asymptotic of $\MP\{\sup_{t \ge 0}(X(t) - ct^{\be} - \zeta) > u\}$ as $u \to \infty$, where $\zeta$ is a random variable, $X, \eta, \zeta$ are independent.}  

\bigskip

\begin{center}

{\bf \Large Section 1. Tail Asymptotics of Sum\\ and Product of Random Variables}

\end{center}

\bigskip

{\bf\Large 1. Introduction.}

\bigskip

All random variables and processes in this article are real-valued.
Recall some well-known basic definitions.

\bigskip

{\bf Definition 1.} {\it The distribution function $F_X$} of a random variable $X$ is a function 
$F_X : \Mr \to \Mr$, $F_X(u) := \MP\{X \le u\}$. {\it The tail} or {\it the survival function} $\SF_X$ 
of a random variable $X$ is a function $\SF_X : \Mr \to \Mr$, $\SF_X(u) := 1 - F_X(u) = \MP\{X > u\}$.

\bigskip

{\bf Definition 2.} {\it The essential supremum} of $X$ (denoted by $\ess\sup X$) is a real number or $+\infty$
defined as $\ess\sup X = \min\{C \in \Mr\mid X \le C\ \mbox{a.s.}\}$, if the set of these $C$ is nonempty, è $+\infty$, 
if it is empty. Sometimes it is denoted by $\vraimax X$. Similarly, {\it the essential infimum} of 
$X$ is a real number or $-\infty$, defined as $\ess\inf X := \max\{C \in \Mr\mid X \ge C\ \mbox{a.s.}\}$.

\bigskip

What is the asymptotic of the tail $\SF_X(u)$ as $u \uparrow \ess\sup X$ (we call it just 
{\it tail asymptotics of} $X$)? This is the classical problem in Probability Theory.

In Section 1, we consider two independent random variables $X$ and $Y$ with the given tail asymptotic. 
What is the tail asymptotic for $X + Y$ and $XY$? Our main tool is the Laplace method used to find the asymptotic of
integral 
$$
\IL_a^bf(x)e^{\lambda S(x)}dx,\ \ \lambda \to +\infty.
$$
In Section 2, we apply these results to find asymptotics of the excursion probability of a given level by a 
conditionally Gaussian process. The main idea in this section to use the self-similarity (which is imposed as an additional
condition).

We consider all random variables and processes on some fixed probability space $(\Oa, \CF, \MP)$. 
All asymptotic and limit relations hold true as $u \to \infty$, unless otherwise stated.
In this article, $\Gamma(\al)$ for $\al > 0$ denotes the Euler gamma function $\Gamma(\al) = \int_0^{+\infty}x^{\al-1}e^{-x}dx$. As usual, $\Mr_+ := [0, \infty)$.

\bigskip

{\bf\Large 2. Tail asymptotic of $X + Y$ for $\ess\sup X = \infty$, $\ess\sup Y = \sigma \in \Mr$}

\bigskip

Let us first consider the case $\ess\sup X = \infty$, $\ess\sup Y = \sigma \in \Mr$. 
Suppose the tail asymptotic of $X$ and $Y$ is given:
$$
\SF_X(u) \backsim h_X(u) := C_Xu^{\ga}\exp\left(-K_Xu^{\al}\right),
$$
and 
$$
\SF_Y(u) \backsim h_Y(u) := C_Y(\sigma - u)^{\mu}, \ \  u \uparrow \sigma.
$$
Here $C_X, C_Y, K_X, \al, \mu > 0$, $\ga \in \Mr$ are constants. Recall $X$ and $Y$ are independent, hence 
$\ess\sup(X + Y) = \infty$.

\bigskip

{\bf Theorem 1.} Suppose $\al > 1$. Then
$$
\SF_{X+Y}(u) \backsim h_{X+Y}(u) := C_{X+Y}u^{\mu + \ga - \al\mu}\exp\left(-K_X(u - \sigma)^{\al}\right),
$$
where for the sake of brevity
$$
C_{X+Y} := C_XC_Y(K_X\al)^{-\mu}\Ga(\mu + 1).
$$

\bigskip

{\it Proof of Theorem 1.} For the sake of simplicity, let $\sigma = 0$. The general case is easily reduced to this particular one. For $u \in \Mr$
$$
\SF_{X+Y}(u) = \MP\{X + Y > u\} = \MP\{X > u - Y\} = \ME\SF_X(u - Y),
$$
since $X$ and $Y$ are independent. Choose $\delta > 0$ (later we shall define the particular value of $\delta$). 
We have:
$$
\ME\SF_X(u - Y) = \ME\SF_X(u - Y)I_{\{Y < -\delta\}} + \ME\SF_X(u - Y)I_{\{-\delta \le Y \le 0\}}.
$$
But $0 \le \ME\SF_X(u - Y)I_{\{Y < -\delta\}} \le \ME\SF_X(u + \delta)I_{\{Y < -\delta\}}$, since for
$Y < -\delta$ we have $u - Y \ge u + \delta$, and $\SF_X$ is nonincreasing. And $\ME\SF_X(u + \delta)I_{\{Y < -\delta\}} \le \SF_X(u + \delta) \backsim h_X(u + \delta) = o(h_{X+Y}(u))$. (This last relation is straightforward to check.)
Therefore, $\ME\SF_X(u - Y)I_{\{Y < -\delta\}} = o(h_{X+Y}(u))$.

It suffices to prove that $\ME\SF_X(u - Y)I_{\{-\delta \le Y \le 0\}} \backsim h_{X+Y}(u)$.

We have: $u - y \to \infty$ and $\SF_X(u - y) \backsim h_X(u - y)$ uniformly for $y \in [-\delta, 0]$.
Using Lemma 11, we obtain:
$$
\ME\SF_X(u - Y)I_{\{-\delta \le Y \le 0\}} \backsim \ME h_X(u - Y)I_{\{-\delta \le Y \le 0\}}.
$$
Let us rewrite $\ME h_X(u - Y)I_{\{-\delta \le Y \le 0\}}$ as a Stieltjes integral:
$$
\ME h_X(u - Y)I_{\{-\delta \le Y \le 0\}} = \IL_{-\delta}^0h_X(u - y)dF_Y(y) = -\IL_{-\delta}^0h_X(u - y)d\SF_Y(y).
$$
(We use the fact that $\SF_Y \equiv 1 - F_Y$.) Integrating by parts, we obtain:
$$
\ME h_X(u - Y)I_{\{-\delta \le Y \le 0\}} = -\biggl[\biggl.h_X(u - y)\SF_Y(y)\biggr|_{y = -\delta}^{y = 0} - \IL_{-\delta}^0\SF_Y(y)dh_X(u - y)\biggr].
$$
Note that boundary terms do not contribute to the asymptotic. Indeed, $\SF_Y(0) = 0$, since $\ess\sup Y = 0$. 
And $h_X(u + \delta)\SF_Y(-\delta) = o(h_{X+Y}(u))$. Thus,
$$
\ME h_X(u - Y)I_{\{-\delta \le Y \le 0\}} = \IL_{-\delta}^0\SF_Y(y)dh_X(u - y) + o(h_{X+Y}(u)).
$$
Take arbitrarily small $\eps > 0$ and find $\delta > 0$ such that for $y \in [-\delta, 0]$ we have
$$
(1 - \eps)C_Y(-y)^{\mu} \le \SF_Y(y) \le (1 + \eps)C_Y(-y)^{\mu}.
$$
Indeed, $\SF_Y(y) \backsim C_Y(-y)^{\mu}$ as $y \uparrow \sigma = 0$.
The function $h_X$ strictly increases for $u \in [u_0, +\infty)$, where $u_0 := (\ga/K_X\al)^{1/\al}$, since
$h'_X(u) = -C_Xu^{\ga - 1}(K_X\al u^{\al} - \ga)\exp\left(-K_Xu^{\al}\right) < 0$ for $u > u_0$. Therefore, the function $y \mapsto h_X(u - y)$ strictly increases on $[-\delta, 0]$ 
if $u > u_0 + \delta$. Hence for $u > u_0 + \delta$
$$
(1 - \eps)C_Y\IL_{-\delta}^0(-y)^{\mu}dh_X(u - y) \le \IL_{-\delta}^0\SF_Y(y)dh_X(u - y) \le (1 + \eps)C_Y\IL_{-\delta}^0(-y)^{\mu}dh_X(u - y).
$$
So we have eliminated the functions $\SF_X, \SF_Y$. We have proved: for any $\eps > 0$, there exists
 $\delta > 0$ such that  for $u > u_0 + \delta$ we have
$$
(1 - \eps)K(u) + o(h_{X+Y}(u)) \le (1 + o(1))\SF_{X+Y}(u) \le (1 + \eps)K(u) + o(h_{X+Y}(u)), \eqno (1)
$$
where
$$
K(u) := C_Y\IL_{-\delta}^0(-y)^{\mu}dh_X(u - y).
$$
Let us find the asymptotic of this integral. For any $u > u_0 + \delta$, the function $y \mapsto h_X(u - y)$ 
is continuously differentiable on $[-\delta, 0]$, and 
$$
\frac{\pa h_X(u - y)}{\pa y} = C_X(u - y)^{\ga - 1}(K_X\al(u - y)^{\al} - \ga)\exp\left(-K_X(u - y)^{\al}\right).
$$
Hence we can write this Riemann-Stieltjes integral as an ordinary Riemann integral:
$$
K(u) = C_Y\IL_{-\delta}^0(-y)^{\mu}\frac{\pa h_X(u - y)}{\pa y}dy = C_XC_YK_X\al I(u; \al, \al + \ga - 1, \mu) - C_XC_Y\ga I(u; \al, \ga - 1, \mu),
$$
where for $\al > 0, \be \in \Mr, \mu > 0, u > 0$ we denote
$$
I(u; \al, \be, \mu) := \IL_{-\delta}^0(-y)^{\mu}(u - y)^{\be}\exp\left(-K_X(u - y)^{\al}\right)dy =
\IL_{0}^{\delta}z^{\mu}(u + z)^{\be}\exp\left(-K_X(u + z)^{\al}\right)dz.
$$

It suffices to find asymptotics of $I(u; \al, \be, \mu)$. Denote $K := K_X$ for the sake of brevity.

\bigskip

{\bf Lemma 1.}
$$
I(u; \al, \be, \mu) \backsim g(u) := (K\al)^{-\mu-1}\Ga(\mu + 1)u^{\be - (\al - 1)(\mu + 1)}\exp\left(-Ku^{\al}\right).
$$

\bigskip

{\it Proof of Lemma 1.} One cannot directly apply the Laplace method (see [2]), since 
the exponent does not contain $u$ as a multiplier, it depends on $u$ in a more complex way.
Let us get rid of the multiplier $(u + z)^{\be}$ in the integrand. 
Notice that $(u + z)^{\be} \backsim u^{\be}$ uniformly for $z \in [0, \delta]$. By Lemma 11
$$
I(u; \al, \be, \mu) \backsim u^{\be}I(u; \al, \mu),
$$
where
$$
I(u; \al, \mu) := I(u; \al, 0, \mu).
$$
Change the variable in $I(u; \al, \mu)$ to eliminate the cumbersome exponent,
but to "preserve the scale" of its dependence on $u$. 
$$
w = \left((u + z)^{\al} - u^{\al}\right)/(\al u^{\al - 1}),\ z = (\al u^{\al - 1}w + u^{\al})^{1/\al} - u.
$$
Then $dz = u^{\al - 1}(\al u^{\al - 1}w + u^{\al})^{1/\al - 1}dw$. In particular,
$\left.dz\right|_{z=0} = \left.dw\right|_{w=0}$. (This is what we call "preserving the scale".) 
The integration segment $[0, \delta]$ maps into $[0, \bd(u)]$, 
$\bd(u) := ((u + \delta)^{\al} - u^{\al})/(\al u^{\al - 1}) \to \delta$. Therefore, we have:
$$
I(u; \al, \mu) = \IL_0^{\bd(u)}\left[\left(\al u^{\al - 1}w + u^{\al}\right)^{1/\al} - u\right]^{\mu}u^{\al - 1}
\exp\left(-K\al u^{\al-1}w - Ku^{\al}\right)\left(\al u^{\al-1}w + u^{\al}\right)^{1/\al - 1}dw =
$$
$$
= u^{\al - 1}\exp\left(-Ku^{\al}\right)\IL_0^{\bd(u)}\left[\left(\al u^{\al - 1}w + u^{\al}\right)^{1/\al} - u\right]^{\mu}\exp\left(-K\al u^{\al-1}w\right)\left(\al u^{\al-1}w + u^{\al}\right)^{1/\al - 1}dw =
$$
$$
= u^{\al - 1}\exp\left(-Ku^{\al}\right)u^{\mu(\al - 1)/\al}u^{(\al - 1)(1/\al - 1)}\IL_0^{\bd(u)}\left[(\al w + u)^{1/\al} - u^{1/\al}\right]\exp\left(-K\al u^{\al-1}w\right)(\al w + u)^{1/\al - 1}dw =
$$
$$
= \exp\left(-Ku^{\al}\right)u^{(\al - 1)(\mu + 1)/\al}\al^{-1}\IL_0^{\al\bd(u)}\left[(v + u)^{1/\al} - u^{1/\al}\right]^{\mu}\exp\left(-Ku^{\al-1}v\right)(v + u)^{1/\al - 1}dv.
$$
(We have changed variables again: $v = \al w$.) $(v + u)^{1/\al - 1} \backsim u^{1/\al - 1}$ uniformly for $v \ge 0$.
Hence by Lemma 11 we can substitute $(v + u)^{1/\al - 1}$ by $u^{1/\al - 1}$ in the integrand, 
and this will not change the asymptotic. Thus
$$
I(u; \al, \mu) \backsim \al^{-1}\exp\left(-Ku^{\al}\right)u^{(\al - 1)\mu/\al}J_{\al\bd(u)}(u; \al, \mu).
$$
Here for $\delta' > 0$ 
$$
J_{\delta'}(u; \al, \mu) := \IL_0^{\delta'}\left[(v + u)^{1/\al} - u^{1/\al}\right]^{\mu}\exp\left(-Ku^{\al-1}v\right)dv.
$$
Let us find the asymptotic of this integral. We shall show that it is the same for all $\delta' > 0$. 
But  $\bd(u) \to \delta$, hence for sufficiently large $u$ $\al\delta/2 < \al\bd(u) < 2\al\delta$, and
$$
J_{\al\delta/2}(u) \le J_{\al\bd(u)}(u) \le J_{2\al\delta}(u),
$$
and $J_{\al\bd(u)}(u)$ has the same asymptotic as $J_{\delta'}(u),\ \delta' > 0$.

By Lemma 11, we can replace  $\left[(v + u)^{1/\al} - u^{1/\al}\right]^{\mu}$ by $(\al^{-1}u^{1/\al - 1}v)^{\mu}$ in the integrand of $J_{\delta'}(u; \al, \mu)$, since
$$
[(v + u)^{1/\al} - u^{1/\al}]^{\mu} \backsim [\al^{-1}u^{1/\al - 1}v]^{\mu}
$$
uniformly for $v \in [0, \delta']$. Let us prove this asymptotic relation. Use the Taylor expansion for the function 
$u \mapsto u^{1/\al}$
$$
(v + u)^{1/\al} - u^{1/\al} = \al^{-1}u^{1/\al - 1}v + \al^{-1}(\al^{-1} - 1)(u + \theta v)^{1/\al - 2}v^2,
$$
where $\theta \in [0, 1]$ depends on $u, v$. Since $\al > 1$, we have $1/\al - 2 < 0$ and $(u + \theta v)^{1/\al - 2} \le 
u^{1/\al - 2}$. Therefore,
$$
|(u + \theta v)^{1/\al - 2}v^2| \le u^{1/\al - 2}\delta'v = o\left(\al^{-1}u^{1/\al - 1}v\right)
$$
uniformly for $v \in [0, \delta']$. Hence $(v + u)^{1/\al} - u^{1/\al} \backsim \al^{-1}u^{1/\al - 1}v$.

Therefore,
$$
J_{\delta'}(u; \al, \mu) \backsim \al^{-\mu}u^{(1/\al - 1)\mu}\IL_0^{\delta'}v^{\mu}\exp\left(-Ku^{\al-1}v\right)dv.
$$
We have $u^{\al - 1} \to \infty$ for $\al > 1$. By Watson's lemma (see [2])
$$
\IL_0^{\delta'}v^{\mu}\exp\left(-uv\right)dv \backsim u^{-\mu-1}\Ga(\mu + 1).
$$
Thus
$$
J_{\delta'}(u; \al, \mu) \backsim \al^{-\mu}\Ga(\mu + 1)u^{(1/\al - 1)\mu}(Ku^{\al - 1})^{-\mu-1} = \al^{-\mu}K^{-\mu-1}\Ga(\mu + 1)u^{(1/\al - 1)\mu - (\al - 1)(\mu + 1)}.
$$
Recall that $J_{\al\bd(u)}(u)$ has the same asymptotic. After easy technical calculations we obtain:
$$
I(u; \al, \mu) \backsim (K\al)^{-\mu-1}\Ga(\mu+1)u^{-(\al - 1)(\mu + 1)}\exp\left(-Ku^{\al}\right),
$$
and
$$
I(u; \al, \be, \mu) \backsim (K\al)^{-\mu-1}\Ga(\mu+1)u^{\be - (\al - 1)(\mu + 1)}\exp\left(-Ku^{\al}\right). \ \ \square
$$

\bigskip

{\it Proof of Theorem 1.}  From Lemma 1, we immediately obtain: $I(u; \al, \ga - 1, \mu) = o(I(u; \al, \al + \ga - 1, \mu))$.
Hence
$$
K(u) \backsim C_XC_YK_X\al I(u; \al, \al + \ga - 1, \mu) \backsim
$$
$$
\backsim C_XC_YK_X\al (K_X\al)^{-\mu-1}\Ga(\mu+1)u^{\al + \ga - 1 - (\al - 1)(\mu + 1)}\exp\left(-Ku^{\al}\right) = h_{X+Y}(u).
$$
Dividing (1) by $h_{X+Y}(u)$, we get:
$$
1 - \eps \le \varliminf\limits_{u\to\infty}\frac{\SF_{X+Y}(u)}{h_{X+Y}(u)}  \le \varlimsup\limits_{u\to\infty}\frac{\SF_{X+Y}(u)}{h_{X+Y}(u)} \le 1 + \eps.
$$
Since $\eps > 0$ is arbitrary,
$$
\lim\limits_{u\to\infty}\frac{\SF_{X+Y}(u)}{h_{X+Y}(u)} = 1.
$$
The case $\sigma = 0$ is proved. The general case is reduced to this one by the 
obvious change of variables: $\tilde Y := Y - \sigma, \tilde u := u - \sigma$. $\tilde u \to \infty$ 
as $u \to \infty$, hence
$$
\SF_{X+Y}(u) = \MP\{X + Y > u\} = \MP\{X + \tilde Y > \tilde u\} \backsim h_{X + \tilde Y}(\tilde u),
$$
because $\ess\sup\tilde Y = 0$. It can be easily shown that
$$
h_{X + \tilde Y}(\tilde u) \backsim h_{X+Y}(u),
$$
since $\tilde u^{\mu + \ga - \al\mu} = (u - \sigma)^{\mu + \ga - \al\mu} \backsim u^{\mu + \ga - \al\mu}$. 
The proof is complete. $\square$

\bigskip

{\bf\Large 3. Tail asymtotic of $X + Y$ for $\ess\sup X = \ess\sup Y = \infty$}

\bigskip

Let $\ess\sup X = \ess\sup Y = \infty$. We do not need to specify any particlular type of asymptotic for $X$, $Y$;
the results of this subsection are valid for a fairly broad class of asymptotic.
Let us introduce some additional conditions. 

\bigskip

{\bf Definition 3.} Denote by $\CM$ the class of all functions $f : \Mr_+ \to \Mr_+$ with the two following properties:

1. there exists $u_0 > 0$ such that $f$ is nonincreasing on $[u_0, \infty)$;

2. $f(u) \to 0$.

\bigskip

{\bf Remark 1.} The survival function of any random variable is in $\CM$.

\bigskip

{\bf Definition 4.} Suppose $f, g \in \CM$. Then:

- the ordered pair $(f, g)$ {\it satisfies the} (A) {\it condition} if there exists a function 
$\qi : \Mr_+ \to \Mr_+$ such that $\qi(u) \to \infty$, $\qi(u)/u \to 0$, $f(\qi(u)) = o(g(u))$ 
and $g(u) = g(u - \qi(u))$.

- the ordered pair $(f, g)$ {\it satisfies the} (B) {\it condition} if there exists a function 
$\qi : \Mr_+ \to \Mr_+$ such that  $\qi(u) \to \infty$, $\qi(u)/u \to 0$, $f(\qi(u)) = o(g(u))$
and $g(u - \qi(u)) = g(u + \qi(u))$.

\bigskip

{\bf Remark 2.} If $f_1, f_2, g_1, g_2 \in \CM$ and $f_1 \backsim f_2$, $g_1 \backsim g_2$, 
then $(f_1, g_1)$ and $(f_2, g_2)$ either both satisfy or both do not satisfy the (A) condition, 
and either both satisfy or both do not satisfy the (B) condition

\bigskip

{\bf Remark 3.} Each of these conditions implies $f(u) = o(g(u))$, since for $u$ large enough 
$\qi(u) < u$, and $f$ is nonincreasing on $[u_0, \infty)$ for sufficiently large $u_0$.

\bigskip

{\bf Remark 4.} Suppose $f_1, f_2, g \in \CM$, $f_1 = o(f_2)$.
If $(f_2, g)$ satisfies any of the conditions (A), (B), then $(f_1, g)$ satisfies it.

\bigskip

{\bf Theorem 2.} Suppose one of the following conditions holds:

1. The pair $(\tF_X, \SF_Y)$ satisfies (B), where $\tF_X : \Mr_+ \to \Mr$ is defined as follows: 
$\tF_X(u) := \SF_X(u) + F_X(-u)$.

2. $X \ge 0$ a.s., and the pair $(\SF_X, \SF_Y)$ satisfies (A).

Then
$$
\SF_{X+Y}(u) \backsim \SF_Y(u).
$$

{\bf Proof of Theorem 2.} Let us prove that the second condition is sufficient. 
Since $X$ and $Y$ are independent,
$$
\SF_{X+Y}(u) = \MP\{X+Y > u\} = \MP\{Y > u - X\} = \ME\SF_Y(u - X).
$$
But
$$
\ME\SF_Y(u - X) = \ME\SF_Y(u - X)I_{\{0 \le X \le \qi(u)\}} + \ME\SF_Y(u - X)I_{\{X > \qi(u)\}}.
$$
The second summand is between $0$ and $\ME I_{\{X > \qi(u)\}}$, because $0 \le \SF_Y(y) \le 1$ for $y \in \Mr$. 
Since $\ME I_{\{X > \qi(u)\}} = \SF_X(\qi(u)) = o(\SF_Y(u))$, the second summand is also $o(\SF_Y(u))$. 
And the first summand is between $\SF_Y(u - \qi(u))$ and $\SF_Y(u)$ since $\SF_Y$ is nonincreasing. 
It suffices to note that $\SF_Y(u - \qi(u)) \backsim \SF_Y(u)$. This completes the proof of the second statement.
The first one is proved similarly, we need to decompose
$$
\ME\SF_Y(u - X) = \ME\SF_Y(u - X)I_{\{|X| \le \qi(u)\}} + \ME\SF_Y(u - X)I_{\{|X| > \qi(u)\}}. \ \ \  \square
$$
\bigskip

How to apply this theorem? Which functions $f, g$ satisfy these conditions?

\bigskip

{\bf Lemma 2.} 1. Let $f(u) := C_1u^{-\al_1},\ g(u) := C_2u^{-\al_2}$, where $C_1, C_2, \al_1, \al_2 > 0$ 
are constants. Then ($(f, g)$ satisfies (A)) $\Leftrightarrow$ ($(f, g)$ satisfies (B)) $\Leftrightarrow$
$\al_1 > \al_2$.

2. Suppose $f(u) := C_1u^{\ga}\exp\left(-Ku^{\al}\right),\ g(u) := C_2u^{-\mu}$, where 
$C_1, C_2, K, \al, \mu > 0$, $\ga \in \Mr$ are constants. Then $(f, g)$ satisfies (A) and (B).

\bigskip

{\bf Proof of Lemma 2.} The second statement immediately follows from the first (see Remark 3).
But if $\al_1 \le \al_2$ then $f(u) \ne o(g(u))$ and neither (A) nor (B) holds true (see Remark 2). 
For $\al_1 > \al_2$, take $\qi(u) := u^{(\al_1 + \al_2)/(2\al_1)}$. $\square$

\bigskip

{\bf\Large 4. Tail asymptotic of $XY$ for $\ess\sup X = \infty$, $\ess\sup Y = \sigma$, $X, Y > 0$}

\bigskip

For the sake of simplicity suppose $X, Y > 0$ a.s. Suppose $\ess\sup X = \infty,\ \ess\sup Y = \sigma$, as in subsection 2. Naturally, $\sigma > 0$. 
Suppose we are given tail asymptotic of $X$ and $Y$, the same as in subsection 2:
$$
\SF_X(u) \backsim h_X(u) := C_Xu^{\ga}\exp\left(-K_Xu^{\al}\right),
$$
and
$$
\SF_Y(u) \backsim h_Y(u) := C_Y(\sigma - u)^{\mu},\ \ u \uparrow \sigma
$$
Here $C_X, C_Y, K_X, \al, \mu > 0$, $\ga \in \Mr$ are constants. Recall $X$ and $Y$ are independent. 
Therefore, $\ess\sup(XY) = \infty$. We need not impose the condition $\al > 1$. 

\bigskip

{\bf Theorem 3.} Under these conditions,
$$
\SF_{XY}(u) \backsim h_{XY}(u) := C_{XY}u^{\ga - \al\mu}\exp\left(-K_X\sigma^{-\al}u^{\al}\right).
$$
Here
$$
C_{XY} := C_XC_Y\Ga(\mu + 1)\sigma^{\al\mu + \mu - \ga}(K_X\al)^{-\mu}.
$$

\bigskip

{\bf Proof of Theorem 3.} The proof is simpler than in subsection 2. Everywhere in this proof $u > 0$. We have
$$
\SF_{XY}(u) = \MP\{XY > u\} = \MP\{X > u/Y\} = \ME\SF_X(u/Y),
$$
since $X$, $Y$ are independent. Fix $\delta \in (0, \sigma)$ (we shall determine its exact value later). 
We have:
$$
\ME\SF_X(u/Y) = \ME\SF_X(u/Y)I_{\{\sigma - \delta \le Y \le \sigma\}} + \ME\SF_X(u/Y)I_{\{Y < \sigma - \delta\}}.
$$
The function $\SF_X$ is nonincreasing, and $y \mapsto \SF_X(u/y)$ is nondecreasing for $y > 0$. Therefore,
$$
0 \le \ME\SF_X(u/Y)I_{\{Y < \sigma - \delta\}} \le \ME\SF_X(u/(\sigma - \delta))I_{\{Y < \sigma - \delta\}} \le \SF_X(u/(\sigma - \delta)) \backsim h_X(u/(\sigma - \delta)) = o(h_{XY}(u)).
$$
(It is easy to verify the last relation.) Hence $\ME\SF_X(u/Y)I_{\{Y < \sigma - \delta\}} = o(h_{XY}(u))$.
It suffices to prove: $\ME\SF_X(u/Y)I_{\{\sigma - \delta \le Y \le \sigma\}} \backsim h_{XY}(u)$.
We have: $u/y \to \infty$, $\SF_X(u/y) \backsim h_X(u/y)$ uniformly for $y \in [\sigma - \delta, \sigma]$. By Lemma 11,
$$
\ME\SF_X(u/Y)I_{\{\sigma - \delta \le Y \le \sigma\}} \backsim \ME h_X(u/Y)I_{\{\sigma - \delta \le Y \le \sigma\}}.
$$
Rewrite this as a Stieltjes integral:
$$
\ME h_X(u/Y)I_{\{\sigma - \delta \le Y \le \sigma\}} = \IL_{\sigma - \delta}^{\sigma}h_X(u/y)dF_Y(y) = - \IL_{\sigma - \delta}^{\sigma}h_X(u/y)d\SF_Y(y) =
$$
$$
= -\left[h_X(u/y)\SF_Y(y)\biggl.\biggr|_{y = \sigma - \delta}^{y = \sigma} - \IL_{\sigma - \delta}^{\sigma}\SF_Y(y)dh_X(u/y)\right].
$$
(We integrated this Stieltjes integral by parts.) The boundary terms do not contribute to the asymptotic, since they are
 $o(h_{XY}(u))$. Indeed, $\ess\sup Y = \sigma$, $\SF_Y(\sigma) = 0$; and 
$0 \le h_X(u/(\sigma - \delta))\SF_Y(\sigma-\delta) \le h_X(u/(\sigma - \delta)) = o(h_{XY}(u))$.

Hence
$$
\SF_{XY}(u) = o(h_{XY}(u)) + \IL_{\sigma - \delta}^{\sigma}\SF_Y(y)dh_X(u/y).
$$
But $h_X \in C^1(0, +\infty)$, $h_X$ is nonincreasing on $[u_0, \infty)$, $u_0 := (\ga/K_X\al)^{1/\al}$. 
(Recall the proof of Theorem 1.) Hence the function $y \mapsto h_X(u/y)$ is nondecreasing on $[\sigma - \delta, \sigma]$ if $u > u_0\sigma$. Now let us define $\delta$. Take an arbitrary $\eps > 0$ 
and find $\delta > 0$ such that $y \in [\sigma - \delta, \sigma]$
$$
(1 - \eps)C_Y(\sigma - y)^{\mu} \le \SF_Y(y) \le (1 + \eps)C_Y(\sigma - y)^{\mu}.
$$
Then for $u > u_0\sigma$
$$
(1 - \eps)I(u) \le \IL_{\sigma - \delta}^{\sigma}\SF_Y(y)dh_X(u/y) \le (1 + \eps)I(u),
$$
where
$$
I(u) := C_Y\IL_{\sigma - \delta}^{\sigma}(\sigma - y)^{\mu}dh_X(u/y).
$$
Thus: for $u > u_0\sigma$
$$
(1 - \eps)I(u) + o(h_{XY}(u)) \le (1 + o(1))\SF_{XY}(u) \le (1 + \eps)I(u) + o(h_{XY}(u)). \eqno(2)
$$
It suffices to find the asymptotic of $I(u)$. We calculated 
$h'_X(u) = -C_Xu^{\ga - 1}(K_X\al u^{\al} - \ga)\exp\left(-K_Xu^{\al}\right)$ during the proof of Theorem 1. 
Hence the function $y \mapsto h_X(u/y)$ is continuously differentiable on $[\sigma - \delta, \sigma]$ and
$$
\frac{\pa h_X(u/y)}{\pa y} = -(u/y^2)\left.h'_X(x)\right|_{x = u/y} = C_X(u/y^2)(u/y)^{\ga - 1}(K_X\al(u/y)^{\al} - \ga)\exp\left(-K_Xy^{-\al}u^{\al}\right).
$$
The Riemann-Stieltjes integral can be rewritten as a Riemann integral:
$$
I(u) := C_XC_Y\IL_{\sigma - \delta}^{\sigma}u^{\ga}y^{-\ga-1}(\sigma - y)^{\mu}\left(K_X\al y^{-\al}u^{\al} - \ga\right)\exp\left(-K_Xy^{-\al}u^{\al}\right)dy =
$$
$$
= C_XC_YK_X\al u^{\al + \ga}I(u^{\al}; \al, -\ga-\al-1, \mu) - C_XC_Y\ga u^{\ga}I(u^{\al}; \al, -\ga - 1, \mu), \eqno (3)
$$
where for $\al > 0, \be \in \Mr, \mu > 0, u > 0$ we denote
$$
I(u; \al, \be, \mu) := \IL_{\sigma - \delta}^{\sigma}y^{\be}(\sigma - y)^{\mu}\exp\left(-K_Xy^{-\al}u\right)dy
= \IL_0^{\delta}(\sigma - z)^{\be}z^{\mu}\exp\left(-K_X(\sigma - z)^{-\al}u\right)dz.
$$
(We changed variables $z := \sigma - y$.) Applying Lemma 10, one can easily find the asymptotic of this integral:
$$
I(u; \al, \be, \mu) \backsim \sigma^{\be + (\al + 1)(\mu + 1)}(K_X\al)^{-\mu-1}\Ga(\mu + 1)u^{-\mu - 1}\exp\left(-K_X\sigma^{-\al}u\right).
$$
Indeed, the function $S(z) := K_X(\sigma - z)^{-\al}$ is strictly increasing on $[0, \delta]$, 
and $S'(0) = K_X\al\sigma^{-\al-1}$; and it suffices to apply Lemma 10 for this $S$ and $f(z) := (\sigma - z)^{\be}$, 
$\mu := \mu + 1$. The asymptotic of $I(u; \al, \be, \mu)$ depends on $\be$ only by the coefficient. 
Hence the second summand in (3) is infinitesimally small with respect to the first summand.
$$
I(u) \backsim C_XC_YK_X\al u^{\al + \ga}\sigma^{-\ga-\al-1 + (\al + 1)(\mu + 1)}(K_X\al)^{-\mu-1}\Ga(\mu+1)(u^{\al})^{-\mu - 1}\exp\left(-K_X\sigma^{-\al}u^{\al}\right) = h_{XY}(u).
$$
Divide (2) by $h_{XY}(u)$ and obtain:
$$
1 - \eps \le \varliminf\limits_{u \to \infty}\frac{\SF_{XY}(u)}{h_{XY}(u)} \le  \varlimsup\limits_{u \to \infty}\frac{\SF_{XY}(u)}{h_{XY}(u)} \le 1 + \eps.
$$
It suffices to note that $\eps > 0$ is arbitrary. The proof is complete. $\square$

\bigskip

{\bf\Large 5. Tail asymptotic of $XY$ for $\ess\sup X = \ess\sup Y = \infty$, $X, Y > 0$}

\bigskip

Again, suppose $X, Y$ are a.s. strictly positive.
Here we need a power tail asymptotic of one of these variables, e.g. $Y$:
$$
\SF_Y(u) \backsim C_Yu^{-\al},
$$
where $C_Y, \al > 0$ are constants. Also, let $\ME X^{\al} < \infty$. Suppose $\SF_X$ 
satisfies the following condition $(\mathrm C_{\al})$:

\bigskip

{\bf Definition 5.} Let $\al > 0$. A function $f \in \CM$ {\it satisfies the $(\mC_{\al})$ condition} 
if there exists a function $\qi : \Mr_+ \to \Mr_+$ such that
 $\qi(u) \to 0,\ \qi(u)/u \to 0,\ f(\qi(u)) = o(u^{-\al})$.

\bigskip

{\bf Remark 5.} $f \in \CM$ satisfies $(\mC_{\al})$ iff the ordered pair $(f(u), u^{-\al})$ satisfies (A).

\bigskip

We also need the condition $\ME X^{\al} < \infty$. Since we operate with survival functions,
let us rewrite this condition in terms of $\SF_X$.

\bigskip

{\bf Definition 6.} Suppose $\al > 0$. The function $f \in \CM$ {\it satisfies the $(\mD_{\al})$ condition} 
if $f(u) = o(u^{-\al})$ and
$$
\IL_0^{\infty}f(u)u^{\al-1}du < \infty.
$$

{\bf Lemma 4.} Suppose $X$ is a nonnegative random variable. Then $\ME X^{\al} < \infty$ iff the function
$\SF_X$ satisfies $(\mD_{\al})$.

\bigskip

{\it Proof of Lemma 4.} It is well-known from classical probability theory, nevertheless we 
expose it in detail. Suppose $\ME X^{\al} < \infty$. First, let us prove that $\SF_X(u) = o(u^{-\al})$. For $u > 0$
$$
u^{\al}\SF_X(u) = u^{\al}\MP\{X > u\} = \ME u^{\al}I_{\{X > u\}} \le \ME X^{\al}I_{\{X > u\}} \to 0
$$
by the Lebesgue dominated convergence theorem, since $\ME X^{\al} < \infty$. 
Rewrite this expectation as a Riemann-Stieltjes integral:
$$
\IL_0^{\infty}u^{\al}dF_X(u) = -\IL_0^{\infty}u^{\al}d\SF_X(u)
$$
and integrate by parts:
$$
- u^{\al}\SF_X(u)\biggl.\biggr|_{u=0}^{u=\infty} + \IL_0^{\infty}\SF_X(u)du^{\al}.
$$
Boundary terms are zero, since $u^{\al}\SF_X(u)\left.\right|_{u=\infty} = 0$ (we just proved this). Hence:
$$
\IL_0^{\infty}\SF_X(u)u^{\al-1}du = \al^{-1}\IL_0^{\infty}\SF_X(u)du^{\al} < \infty.
$$
Hence $\ME X^{\al} < \infty$ implies that the survival function $\SF_X$ satisfies $(\mD_{\al})$. The proof 
of the converse statement is similar. $\square$

\bigskip

{\bf Remark 6.} (Analogous to Remark 2.) If $f, g \in \CM$, $f \backsim g$, then these functions
either both satisfy or both do not satisfy the condition $(\mathrm C_{\al})$, 
and they both satisfy or do not satisfy the condition $(\mD_{\al})$.

\bigskip

{\bf Remark 7.} (Analogous to Remark 3.) The condition $(\mathrm C_{\al})$ implies $f(u) = o(u^{-\al})$.

\bigskip

{\bf Remark 8.} (Analogous to Remark 4.) For $f, g \in \CM$, $f(u) = o(g(u))$, if $g$ satisfies $(\mathrm C_{\al})$,
then $f$ also satisfies this condition; if $g$ satisfies $(\mathrm D_{\al})$, 
then $f$ also satisfies this condition.

\bigskip

{\bf Theorem 4.} Suppose $\SF_Y(u) \backsim C_Yu^{-\al}$, and $\SF_X$ 
satisfies the conditions $(\mC_{\al})$ and $(\mD_{\al})$. Then
$$
\SF_{XY}(u) \backsim C_Y\ME X^{\al}u^{-\al}.
$$

\bigskip

{\it Proof of Theorem 4.} In this proof, $u > 0$. As before, $\SF_{XY}(u) = \ME\SF_Y(u/X)$ for $u > 0$,
since $X$, $Y$ are independent. We obviously have:
$$
\ME\SF_Y(u/X) = \ME\SF_Y(u/X)I_{\{0 < X \le \qi(u)\}} + \ME\SF_Y(u/X)I_{\{X > \qi(u)\}}.
$$
Since $0 \le \SF_Y(u) \le 1$ for all $u$, we have
$$
0 \le \ME\SF_Y(u/X)I_{\{X > \qi(u)\}} \le \MP\{X > \qi(u)\} = \SF_X(\qi(u)) = o(u^{-\al})
$$
(according to the $(\mathrm C_{\al})$ condition). Hence $\ME\SF_Y(u/X)I_{\{X > \qi(u)\}} = o(u^{-\al})$. Applying Lemma
11 and noting that  $u/x \to \infty$ uniformly for $x \in (0, \qi(u)]$ (since $u/\qi(u) \to \infty$), 
$\SF_Y(u/X) \backsim C_Y(u/X)^{-\al}$, we get:
$$
\ME\SF_Y(u/X)I_{\{0 < X \le \qi(u)\}} \backsim \ME C_Y(u/X)^{-\al}I_{\{0 < X \le \qi(u)\}}.
$$
Thus,
$$
\ME\SF_Y(u/X)I_{\{0 < X \le \qi(u)\}} \backsim C_Yu^{-\al}\ME X^{\al}I_{\{0 < X \le \qi(u)\}}.
$$
But $\qi(u) \to \infty$, $\ME X^{\al} < \infty$, hence $\ME X^{\al}I_{\{0 < X \le \qi(u)\}} \to \ME X^{\al}$ 
by the Lebesgue dominated convergence theorem. Thus, 
$\ME\SF_Y(u/X)I_{\{0 < X \le \qi(u)\}} \backsim C_Y\ME X^{\al}u^{-\al}$, 
and the proof is complete. $\square$

\bigskip

What examples of functions $f \in \CM$ satisfying the conditions $(\mathrm C_{\al})$ and $(\mathrm D_{\al})$ 
are there? 

\bigskip

{\bf Lemma 4.} 1. The function $f(u) := C_fu^{-\be}$, where $C_f, \be > 0$ are constants, satisfies
the conditions $(\mathrm C_{\al})$ and $(\mathrm D_{\al})$ for $\be > \al$ and does not satisfy them for $\be \le \al$.

2. The function $f(u) := C_fu^{\ga}\exp\left(-Ku^{\be}\right)$, where $C_f, K, \be > 0$, $\ga \in \Mr$ are constants,
satisfies the conditions $(\mathrm C_{\al})$ è $(\mathrm D_{\al})$.

\bigskip

{\it Proof of Lemma 4.} Remark 7 shows that the second statement follows from the first one.
Let us prove the first statement. Speaking about the condition $(\mC_{\al})$, it suffices
to use Lemma 2 and Remark 5. The condition $(\mD_{\al})$ is not satisfied for $\be \le \al$ can can be straightforwardly
checked for $\be > \al$. $\square$

\bigskip

\begin{center}

{\bf\Large Section 2. Tail asymptotics of extrema\\ of conditionally Gaussian processes}

\end{center}

\bigskip

{\bf\Large 6. Introduction.}

\bigskip

We shall apply this theory to find asymptotic of
$$
\MP\{\sup\limits_{t \ge 0}(X(t) - \eta t^{\be}) > u\}, \ \MP\{\sup\limits_{t \ge 0}(X(t) - \eta t^{\be} - \zeta) > u\},
$$
where $X = (X(t), t \ge 0)$ is a Gaussian centered self-similar locally stationary process (we shall clarify 
the conditions imposed on $X$ later), $\eta > 0, \zeta$ - independent random variables, $(\eta, \zeta)$ 
is independent of $X$.

What is the history of this problem?

Classical asymptotical theory of extrema of Gaussian processes and fields was developed (by, e.g., Piterbarg
and Pickands, see monograph [3]) for centered processes and fields. But afterwards, non-centered process (i.e.
processes with a trend) were considered. They have the form $(X(t) + m(t))$, where $X$ 
is a centered Gaussian processes, and $m$ is a nonzero deterministic function, which is called a {\it trend}. 
See, e.g., a well-known article [1], where $m(t) = -ct^{\be},\ t \ge 0, \ c, \be$ are constants ({\it power trend}).

The problem was then generalized to the case of {\it conditionally Gaussian processes} $Y$.
They depend on random variables $\eta_1, \ldots, \eta_n$ and on a Gaussian process $X$, where 
$(\eta_1, \ldots, \eta_n)$ is independent of $X$. $Y$ is called so because 
the conditional distribution of $Y$ for fixed $\eta_1, \ldots, \eta_n$ is Gaussian.

For example, the following model is considered in [4]: $Y = (Y(t), t \ge 0)$, $Y(t) = X(t)(\eta - \zeta t^{\al})$,
where $(\eta, \zeta)$ is independent of $X$, the random variables $\eta, \zeta$ are positive, bounded and 
$\ess\inf\eta > \eps$.

We shall consider a process $Y = (Y(t), t \ge 0), Y(t) = X(t) - \eta t^{\be}$, where 
$(X(t), t \ge 0)$ is a Gaussian centered process, $\eta > 0$ is independent of $X$, 
but the conditions imposed on $\eta$ are not as strict as the conditions on $\eta, \zeta$ in [4]. 
This model is similar to the one from [1], and we shall intensively use the results from [1]. 
But there is a significant difference: instead of the deterministic trend $-ct^{\beta}$,
we have a random process $(-\eta t^{\beta}, t \ge 0)$. Let us call it a {\it random trend}.

Also, we shall find the asymptotics of
$$
\sup\limits_{t \ge 0}(X(t) - \eta t^{\be} - \zeta),
$$
where $\zeta$ is a random variable, $X, \eta, \zeta$ are independent.

\bigskip

{\bf\Large 7. Basic definitions.}

\bigskip

{\bf Definition 7.} A square-integrable process $X = (X(t), t > 0)$ with $\ME X(t) = 0,\ \ME X^2(t) = 1$ 
is called {\it locally stationary at the point $s > 0$ with the local stationarity index $\al \in (0, 2]$ 
and the limit constant $D(s) > 0$} if
$$
\lim\limits_{t, t' \to s}\frac{\ME(X(t) - X(t'))^2}{|t - t'|^{\al}} = D(s).
$$

\bigskip

{\bf Definition 8.} A random process $X = (X(t), t \ge 0)$ is called {\it self-similar with self-similarity} 
({\it Hurst}) {\it parameter} $H \in (0; 1]$ if for any $a > 0$
$$
(X(at), t \ge 0) \stackrel{d}{=} (a^HX(t), t \ge 0).
$$

\bigskip

{\bf Remark 9.} A Gaussian process $X = (X(t), t \ge 0)$ with $\ME X(t) = 0$ 
is self-similar with Hurst parameter $H$ iff $R(t, t') := \ME X(t)X(t')$ is homogeneous of order $2H$,
i.e. for all $t, t' \ge 0, a > 0$ $R(at, at') = a^{2H}R(t, t')$. 
In particular, its variation is $\ME X^2(t) = ct^{2H}$, where $c = \ME X^2(1)$ is independent of $t$.

\bigskip

{\bf Remark 10.} If a process $X = (X(t), t \ge 0)$ ñ $\ME X(t) = 0, \ME X^2(t) = t^{2H}$ 
is self-similar with Hurst parameter $H$, and the standardized process $Y = (Y(t), t > 0), Y(t) = t^{-H}X(t)$ is locally stationary with index $\al$ and limit constant $D(s_0)$ at the point $s_0$ then it is straightforward to prove:
$Y$ is locally stationary at every point $s > 0$ with the same self-similarity index $\al$, but
with the limit constant $D(s) = (s_0/s)^{\al}D(s_0)$.

\bigskip

{\bf Remark 11.} For the process from the previous remark, $H$ and $\al$ are not related. Changing
$X$ by $(X(t^a), t \ge 0)$, we get a different $H$, but the same $\alpha$.

\bigskip

{\bf Definition 9.} Suppose $H \in (0, 1]$. {\it Fractional Brownian motion with parameter $H$} is a Gaussian
process $B_H = (B_H(t), t \ge 0)$ with a.s. continuous trajectories and the following properties:
$\ME B_H(t) = 0$, $\ME B_H(t)B_H(t') = (t^{2H} + t'^{2H} - |t - t'|^{2H})/2$ for every $t, t' \ge 0$.

\bigskip

{\bf Remark 12.} This is a classic example of a self-similar process with Hurst parameter $H$. 
It is shown in [1] that the process $Y = (Y(t), t > 0),\ Y(t) = t^{-H}B_H(t)$ (note that 
$\ME B_H^2(t) = t^{2H}, t \ge 0$) is locally stationary at every point $s > 0$ with self-similarity index 
$\al = 2H$ and $D(s) = s^{-2H}$. (The reader can easily check this fact himself.)

\bigskip

{\bf Remark 13.} $B_{1/2}$ is a standard Brownian motion.

\bigskip

{\bf\Large 8. Results for a deterministic trend.}

\bigskip

Let us expose the core results from the article [1] in detail, since we shall need them.
Let $H \in (0, 1), c > 0, \al \in (0, 2], \be > H$ be constants. Suppose a stochastic process
$X = (X(t), t \ge 0)$ satisfies the following conditions:

(i) it is Gaussian;

(ii) $\ME X(t) = 0,\ \ME X^2(t) = t^{2H}$ for $t \ge 0$;

(iii) $X$ is self-similar with Hurst parameter $H \in (0, 1)$;

(iv) the standardized process $Y = (Y(t), t > 0)$, $Y(t) := t^{-H}X(t)$ is locally stationary at the point $s_0$
with local stationarity index $\al$ and limit constant $D(s_0)$, where
$$
s_0 := \left(\frac H{c(\be - H)}\right)^{1/\be}.
$$

Then, as $u \to \infty$, we have
$$
p(u, c) := \MP\{\sup\limits_{t \ge 0}(X(t) - ct^{\be}) > u\} \backsim f(u, c), \eqno (1)
$$
where for $\al < 2$
$$
f(u, c) := \frac{H_{\al}\sqrt{\pi}(D(s_0))^{1/\al}}{\sqrt{B}2^{1/\al - 1/2}}A^{2/\al - 1/2}u^{(1 - H/\be)(2/\al - 1)}\Psi\left(Au^{1 - H/\be}\right),
$$
and for $\al = 2$
$$
f(u, c) := 2\sqrt{\frac{AD + B}B}\Psi\left(Au^{1 - H/\be}\right).
$$
$\Psi(x)$ is the tail of the standard normal distribution:
$$
\Psi(x) := \frac1{\sqrt{2\pi}}\int\limits_x^{+\infty}e^{-y^2/2}dy.
$$
$H_{\al}$ is a positive constant, called a {\it Pickands constant}:
$$
H_{\al} := \lim\limits_{T \to \infty}\frac1T\ME\exp\max\limits_{0\le t \le T}(\sqrt{2}B_{\al/2}(t) - t^{\al}).
$$
(See [3], section $D$, for the proof that this constant is indeed well-defined and $H_{\al} > 0$.)
Finally, $A$ and $B$ are positive constants:
$$
A := \left(\frac H{c(\be - H)}\right)^{-H/\be}\frac{\be}{\be - H},\ B := \left(\frac H{c(\be - H)}\right)^{-(H+2)/\be} H\be.
$$
We shall rewrite this result in a more convenient way. If $\qi(x) := (2\pi)^{-1/2}e^{-x^2/2}$ is a standard Gaussian density, then it is easy to verify by L'Hospital's rule that $ \Psi(x) \backsim \qi(x)/x$ as $x \to \infty$. (See also, e.g. [8], Section 7.1, Lemma 2.) But $Au^{1 - H/\be} \to \infty$, hence $\Psi(Au^{1 - H/\be}) \backsim 
\qi(Au^{1 - H/\be})/(Au^{1 - H/\be})$. Thus, $f(u, c) \backsim g(u, c)$, where for $\al \in (0, 2)$ we have
$$
g(u, c) := Cu^{(1 - H/\be)(2/\al - 2)}\qi\left(Au^{1 - H/\be}\right),\
C := \frac{H_{\al}\sqrt{\pi}(D(s_0))^{1/\al}}{\sqrt{B}2^{1/\al - 1/2}}A^{2/\al - 3/2}
$$
and for $\al = 2$ we have
$$
g(u, c) := Cu^{H/\be - 1}\qi\left(Au^{1 - H/\be}\right),\
C := 2\sqrt{\frac{AD + B}B}A^{-1}.
$$
To find how $g$ depends on $c$, denote
$$
K_s := \left(\frac H{\be - H}\right)^{1/\be},\ K_A := \left(\frac H{\be - H}\right)^{-H/\be}\frac{\be}{\be - H},\ K_B := \left(\frac H{\be - H}\right)^{-(H+2)/\be} H\be.
$$
Then we get
$$
s_0 = K_sc^{-1/\be},\ A = K_Ac^{H/\be},\ B = K_Bc^{(H+2)/\be}.
$$
The process $Y$ is locally stationary at every point $s > 0$ with the same local stationarity index $\al$ but with
limit constant $D(s) = (s_0/s)^{\al}D(s_0)$. (See Remark 10.) Hence 
$D(s_0) = D(K_s)(c^{-1/\be})^{-\al} = K_Dc^{\al/\be}$, where $K_D := D(K_s)$. Thus
$$
C = Kc^{H/\be(2/\al - 2)},
$$
where for $\al \in (0, 2)$
$$
K := \frac{H_{\al}\sqrt{\pi}K_D^{1/\al}}{K_B^{1/2}2^{1/\al - 1/2}}K_A^{2/\al - 3/2},
$$
and for $\al = 2$
$$
K := \frac2{K_A}\sqrt{\frac{K_AK_D + K_B}{K_B}}.
$$
Finally, we obtain:
$$
g(u, c) = Kc^{H/\be(2/\al - 2)}u^{(1 - H/\be)(2/\al - 2)}\qi\left(K_Ac^{H/\be}u^{1 - H/\be}\right). \eqno (2)
$$

\bigskip

{\bf\Large 9. Asymptotic of $\MP\{\sup\limits_{t \ge 0}(X(t) - \eta t^{\be}) > u\}$.}

\bigskip

Let $\CS_0 := \sup\limits_{t \ge 0}(X(t) - \eta t^{\be})$. Let us find the tail asymptotic of $\CS_0$. 
We need to impose an additional restriction (v) on the process $X$ :

(v) a.s. there exists an $s \ge 0$ such that $X(s) > 0$.

It could possibly be implied by the other conditions (i) - (iv), but we could not prove this.
The standard Brownian motion and the fractional Brownian motion obviously satisfies (v).
(See, e.g., [7], where the law of iterated logarithm for the fractional Brownian motion is proved as Theorem 3.3; 
this law immediately implies the condition (v).)

\bigskip

{\bf Theorem 5.} Suppose a stochastic process $X$ satisfies the conditions (i) - (v). 
Let $\eta > 0$ be a random variable independent of $X$, $\delta = \ess\inf\eta$. Suppose 
$\MP\{\eta < u\} \backsim C_{\eta}(u - \delta)^{\mu}$ as $u \downarrow \delta$, 
where $C_{\eta}, \mu > 0$ are constants. Then we have:

1. for $\delta > 0$:
$$
\SF_{\CS_0}(u) \backsim K_0u^{\nu}\qi(K_A\delta^{H/\be}u^{1 - H/\be}),
$$
where for the sake of brevity
$$
K_0 := \frac{C_{\eta}K\Gamma(\mu + 1)\be^{\mu}}{K_A^{2\mu}H^{\mu}}\delta^{H/\be(2/\al - 2) - \mu(2H/\be - 1)},\ \nu := \left(1 - \frac H{\be}\right)\left(\frac2{\al} - 2 - 2\mu\right);
$$

2. for $\delta = 0$:
$$
\SF_{\CS_0}(u) \backsim C_{\eta}\CE_{H\mu/\be, \be}u^{-\be\mu/H},
$$
where for $\al, \be > 0$ $\CE_{\al, \be} = \CE_{\al, \be}(X)$ is a positive constant:
$$
\CE_{\al, \be}(X) := \ME\left(\sup\limits_{t \ge 0}\frac{X(t)}{1 + t^{\be}}\right)^{\al}.
$$

{\it Proof of Theorem 5.} We follow the proof of the similar theorem in [1]. Let for $s \ge 0, c > 0$
$$
Z_c(s) := \frac{X(s)}{1 + cs^{\be}},\ \tilde Z_c(s) := c^{H/\be}Z_c((2/D)^{1/\al}s).
$$

\bigskip

{\bf Lemma 5.} For all $u > 0$,
$$
p(u, c) = \MP\{\sup\limits_{s \ge 0}Z_c(s) > u^{1 - H/\be}\}.
$$

{\it Proof of Lemma 5.} This is a consequence of self-similarity of $X$ with Hurst parameter $H$. 
We have:
$$
(X(u^{-1/\be}t), t \ge 0) \stackrel{d}{=} (u^{-H/\be}X(t), t \ge 0).
$$
Thus
$$
\MP\{\sup\limits_{s \ge 0}Z_c(s) > u^{1 - H/\be}\} = \MP\{\exists s \ge 0: X(s) > u^{1 - H/\be} + cs^{\be}u^{1 - H/\be}\} =
$$
$$
= \MP\{\exists s \ge 0: X(u^{-1/\be}s) > u^{1 - H/\be} + c(u^{-1/\be}s)^{\be}u^{1 - H/\be}\} =
$$
$$
= \MP\{\exists s \ge 0: u^{-H/\be}X(s) > u^{1 - H/\be} + cs^{\be}u^{- H/\be}\} = \MP\{\exists s \ge 0: X(s) > u + cs^{\be}\} = p(u, c). \square
$$

\bigskip

{\bf Lemma 6.} The distribution of the process $\tilde Z_c = (\tilde Z_c(s), s \ge 0)$ 
does not depend on $c > 0$.

\bigskip

{\it Proof of Lemma 6.} $\tilde{Z}_c$ is a Gaussian process, hence its distribution is uniquely
determined by its mean and covariance functions. For all $c > 0, t \ge 0$ we have $\ME \tilde{Z}_c(t) = 0$ hence
$\ME X(t) = 0$. Hence it suffices to prove that for all $t, t' \ge 0$, $\ME\tilde{Z}_c(t)\tilde{Z}_c(t')$ 
does not depend on $c > 0$. 
$$
\ME \tilde{Z}_c(t)\tilde{Z}_c(t') = \frac{c^{2H/\be}\ME X((2/D)^{1/\al}t)X((2/D)^{1/\al}t')}{(1 + c((2/D)^{1/\al}t)^{\be})(1 + c((2/D)^{1/\al}t')^{\be})}.
$$
But $D = K_Dc^{\al/\be}$ (see subsection 8), $c((2/D)^{1/\al})^{\be} = 2^{\be/\al}cD^{-\be/\al} = 
(2/K_D)^{\be/\al}$ is independent of $c$. Hence the denominator is independent of $c$. 
And the numerator is equal to $c^{2/\be}2(2/D)^{2/\al}\ME X(t)X(t')$ since $X$ is self-similar (see Remark 9). 
But $c^{2/\be}(2/D)^{2/\al} = (2/K_D)^{\al/\be}$ is also independent of $c$. $\square$

\bigskip

{\it Proof of Theorem 5.} Denote $\tilde Z_c$ just as $\tilde Z$. (Only its distribution is important.)
It is clear that
$$
\sup\limits_{s \ge 0}Z_c(s) = c^{-H/\be}\sup\limits_{s \ge 0}\tilde Z(s).
$$
Condition (v) yields that a.s.
$$
\sup\limits_{s \ge 0}Z_1(s) > 0,\ \sup\limits_{s \ge 0}\tilde Z(s) > 0.
$$
Since $X$ and $\eta$ are independent,
$$
\SF_{\CS_0}(u) = \ME p(u, \eta) = \MP\{\eta^{-H/\be}\sup\limits_{s \ge 0}\tilde Z(s) > u^{1 - H/\be}\}.
$$
But $\eta^{-H/\be} > 0$, $\sup\limits_{s \ge 0}\tilde Z(s) > 0$. It suffices to use Theorems 3, 4.
We have:
$$
\MP\{\sup\limits_{t \ge 0}\tilde{Z}(t) > u^{1 - H/\be}\} = p(u, 1) \backsim  g(u, 1) = Ku^{(1 - H/\be)(2/\al - 2)}\qi\left(K_Au^{1 - H/\be}\right).
$$
Therefore,
$$
\MP\{\sup\limits_{t\ge 0}\tilde{Z}(t) > u\} \backsim Ku^{2/\al - 2}\qi(K_Au),
$$
and
$$
\ess\sup(\sup\limits_{t \ge 0}\tilde Z(t)) = \infty.
$$
The random variable $\sup_{t \ge 0}\tilde Z(t)$ plays the role of $X$ from 
Theorems $3$ and $4$. And the random variable $\eta^{-H/\be}$ plays the role of $Y$.
Consider the cases $\delta > 0$ and $\delta = 0$.

\bigskip

Suppose $\delta > 0$. Then $\ess\sup\eta^{-H/\be} = \delta^{-H/\be}$, and $u^{-\be/H} \downarrow \delta$
as $u \uparrow \delta^{-H/\be}$, 
$$
\SF_{\eta^{-H/\be}}(u) = \MP\{\eta^{-H/\be} > u\} = \MP\{\eta < u^{-\be/H}\} \backsim C_{\eta}(\delta - u^{-\be/H})^{\mu} \backsim C_{\eta}(\be H^{-1}\delta^{\be/H + 1})^{\mu}(\delta^{-H/\be} - u)^{\mu},
$$
since the derivative of the function $u \mapsto u^{-\be/H}$ at the point $u = \delta^{-H/\be}$ 
equals $-\be H^{-1}(\delta^{-H/\be})^{-\be/H-1} = -\be H^{-1}\delta^{1 + H/\be}$, and 
$\delta - u^{-\be/H} \backsim \be H^{-1}\delta^{1 + H/\be}(\delta^{-H/\be} - u)$ as $u \downarrow \delta^{-H/\be}$.

In the notation of Theorem 3
$$
\sigma = \delta^{-H/\be},\ C_Y = C_{\eta}(\be H^{-1}\delta^{\be/H + 1})^{\mu},\ C_X := (2\pi)^{-1/2}K, \ga = 2/\al - 2, K_X := K_A^2/2, \al = 2.
$$
It suffices to apply this theorem and simplify the answer.

\bigskip

Suppose $\delta = 0$. Then $\ess\sup\eta^{-H/\be} = \infty$, and
$$
\SF_{\eta^{-H/\be}}(u) = \MP\{\eta^{-H/\be} > u\} = \MP\{\eta < u^{-\be/H}\} \backsim C_{\eta}(u^{-\be/H})^{\mu} = C_{\eta}u^{-\be\mu/H}.
$$
Since (see above)
$$
\SF_{\sup_{t \ge 0}\tilde Z(t)}(u) \backsim Ku^{2/\al - 2}\qi(K_Au) = \frac K{\sqrt{2\pi}}u^{2/\al - 2}\exp\left(-\frac{K_A^2}2u^2\right),
$$
this function satisfies the conditions $(\mathrm C_{\be\mu/H})$ and $(\mathrm D_{\be\mu/H})$ (see Lemma 4). Hence we can apply Theorem 4. The constant
$$
\ME\left(\sup\limits_{t \ge 0}\tilde Z(t)\right)^{\be\mu/H} = \ME\left(\sup\limits_{t \ge 0}Z_1(t)\right)^{\be\mu/H} = \ME\left(\sup\limits_{t \ge 0}\frac{X(t)}{1 + t^{\be}}\right)^{\be\mu/H}
$$
is denoted by $\CE_{\be\mu/H, \be}(X)$. The proof is complete. $\square$

\bigskip

{\bs\Large 10. Asymptotic of $\MP\{\sup\limits_{t \ge 0}(X(t) - \eta t^{\be} - \zeta) > u\}$.}

\bigskip

Now we shall consider 
$$
\CS := \sup\limits_{t \ge 0}(X(t) - \eta t^{\be} - \zeta) = \CS_0 - \zeta,
$$
where $\zeta$ is a random variable, $X, \eta, \zeta$ are independent. Then $\CS_0, \zeta$ are independent. 
We know the tail asymptotic of $\CS_0$ under certain conditions (see the previous subsection). 
And for a given asymptotic $\MP\{\zeta < u\}$ as $u \downarrow \ess\inf\zeta$, 
we know the tail asymptotic of $-\zeta$, and it suffices to apply Theorems 1, 2.
 
\bigskip

{\bf Theorem 6.} Suppose the process $X$ satisfies (i) - (v). Let $\ess\inf\zeta =: \delta_0$.

1.  Suppose $\delta_0 = -\infty$ and $\MP\{\zeta < -u\} \backsim C_{\zeta}u^{-\ga}$, where 
$C_{\zeta} > 0,\ \ga > 0$. If one of the following conditions holds:

(a) $\delta > 0$;

(b) $\delta = 0$ and $\be\mu/H > \ga$,

then
$$
\SF_{\CS}(u) \backsim C_{\zeta}u^{-\ga}.
$$
If we have

(c) $\delta = 0,\ \be\mu/H < \ga$,

then
$$
\SF_{\CS}(u) \backsim \SF_{\CS_0}(u) \backsim C_{\eta}\CE_{\be\mu/H, \be}(X)u^{-\be H/\mu}.
$$

2. Suppose $\delta_0 \in \Mr$ and $\MP\{\zeta < u\} \backsim C_{\zeta}(u - \delta_0)^{\ga}$, 
where $C_{\zeta}, \ga > 0$ are constants. Suppose also that $\delta > 0$, $2H < \be$. Then
$$
\SF_{\CS}(u) \backsim C_{\CS}u^{\nu}\qi\left(K_A\delta^{H/\be}(u + \delta_0)^{1 - H/\be}\right),
$$
where for the sake of brevity
$$
C_{\CS} := C_{\zeta}K_0K_A^{-2\ga}\delta^{-2H\ga/\be}(1 - H/\be)^{-\ga}\Ga(1 + \ga),
$$
$$
\nu := (2/\al - 2 - 2\mu)(1 - H/\be) + \ga - 2(1 - H/\be)\ga.
$$

{\it Proof of Theorem 6.} Apply directly Theorem 1 for the second case, and Theorem 2 for the first case.
In the first case, in (a) and (b)  $\CS_0$ plays the role of $X$, $-\zeta$ plays the role of $Y$. 
And in (c), their roles are reversed:  $\CS_0$ plays the role of $Y$, $-\zeta$ plays the role of $X$. 

In the second case,  $\CS_0$ plays the role of $X$, $-\zeta$ plays the role of $Y$.
The condition $2H < \be$ is necessary to establish the condition $\al > 1$ in Theorem 1. $\square$

\bigskip

{\bf\Large 11. Conclusion.}

\bigskip

The most interesting case is when neither asymptotic of $X$ nor asymptotic of $\eta$ 
dominate. This is probably the toughest case. It is unlikely that two asymptotical expressions 
can be easily combined. Probably the Pickands method of double sums should be applied 
(see monograph [3], section D or chapter 2).

How to calculate $\CE_{\al, \be}$? We can only do this numerically. 
It is unlikely that one can find an exact form for this constant.
We know the exact form only if $X$ is a Brownian motion, $\be = 1$: 

\bigskip

{\bf Lemma 7.} {\it If $X = B$ is a standard Brownian motion, $\be = 1$, then
$$
\CE_{\al, \be} = 2^{-\al/2}\Gamma\left(\frac{\al}2 + 1\right).
$$}

{\it Proof of Lemma 7.} Let
$$
\CX := \max\limits_{s \ge 0}\frac{B(s)}{1 + s}.
$$
Using Lemma 5, we obtain: for all $u > 0$ 
$$
\MP\{\CX > u^{1/2}\} = \MP\{\max\limits_{s \ge 0}(B(s) - s) > u\}
$$
We used Remark 13: for a Brownian motion $H = 1/2$, $1 - H/\be = 1/2$.
$$
\MP\{\max\limits_{s \ge 0}(-B(s) - s) > u\} = \MP\{\max\limits_{s \ge 0}(B(s) - s) > u\} = \MP\{\exists s \ge 0: B(s) - s = u\} = e^{-2u}.
$$
The last equality uses the result for a hitting time of a Brownian motion with a drift from [5] (chapter 3, section 3.5.C, (5.13)). We also use the 
continuity of Brownian paths: if $\max\limits_{s \ge 0}(B(s) - s) > u$, then for some $s \ge 0$ 
$B(s) - s = u$. Hence $\max\limits_{s \ge 0}(B(s) - s)$ has an exponential distribution, and
$$
\MP\{\max\limits_{s \ge 0}(B(s) - s) > u\} = \MP\{\max\limits_{s \ge 0}(B(s) - s) \ge u\} = e^{-2u}.
$$
For $u > 0$ $\MP\{\CX > u\} = e^{-2u^2}$. Of course, for $u \le 0$ we have $\MP\{\CX > u\} = 0$. 
Hence $\CX$ has the Weibull distribution with parameters $(2, 2)$.  Thus (see [6], chapter 21, section 2; of course one can verify it by a simple calculation) $\CE_{\al, 1} = \ME\CX^{\al} = 2^{-\al/2}\Gamma\left(\frac{\al}2 + 1\right)$. $\square$

\bigskip

\begin{center}

{\bf\Large Appendix}

\end{center}

\bigskip

{\bf\Large 12. Other forms of asymptotic conditions on $X$, $Y$, $\eta$, $\zeta$.}

\bigskip

We used the survival functions of $X$, $Y$, $\eta$, $\zeta$ to impose restrictions on them.
But many distributions are defined in terms of density (with respect to the Lebesgue measure).
Can we rewrite these conditions in terms of density?

Yes, we can. The asymptotic conditions on the survival function are more general than the ones on the density.
Hence we can apply any of these Theorems 1-6 if the asymptotic of denstiy is given.

\bigskip

{\bf Lemma 8.} Suppose $X$ is a random variable, $M := \ess\sup X$. Suppose on a certain left neighborhood
$U \subseteq \Mr$ of $M$ the distribution of $X$ has a density $f_X$ (with respect to the Lebesgue measure). 
This means that for any Borel subset $B \subseteq U$ we have:
$$
\MP\{X \in B\} = \IL_Bf_X(u)du.
$$

1. If $M = +\infty$ and $f_X(u) \backsim C_X\al u^{-\al-1}$, where $C_X, \al > 0$ are constants,
then $\SF_X(u) \backsim C_Xu^{-\al}$.

2. If $M = +\infty$ and $f_X(u) \backsim C_Xu^{\be}\exp\left(-K_Xu^{\al}\right)$, where 
$C_X, \al, K_X > 0,\ \be \in \Mr$ are constants, then 
$\SF_X(u) \backsim C_X\al^{-1}K_X^{-1}u^{\be + 1 - \al}\exp\left(-K_Xu^{\al}\right)$.

3. If $M \in \Mr$ and $f_X(u) \backsim C_X\al(M - u)^{\al-1}$, where $C_X, \al > 0$ are constants,
then $\SF_X(u) \backsim C_X(M - u)^{\al}$.

\bigskip

{\it Proof of Lemma 8.} If $f, g : [a, b) \to \Mr$ are Lebesgue integrable on $[a, b) \subset \Mr$, where
$b \in (a, +\infty]$, and $f(x) \backsim g(x), x \uparrow b$, then by L'Hospital's rule
$$
\IL_x^bf(t)dt \backsim \IL_x^bg(t)dt,\ x \uparrow b.
$$
For $u \le M$ sufficiently close to $M$ we have
$$
\SF_X(u) = \IL_u^Mf_X(t)dt,
$$
Hence statements 1 and 3 are obvious. Statement 2: it suffices to prove that
$$
I(u) := \IL_u^{\infty}C_Xt^{\be}\exp\left(-K_Xt^{\al}\right)dt
\backsim C_X\al^{-1}K_X^{-1}u^{\be + 1 - \al}\exp\left(-K_Xu^{\al}\right). \eqno (4)
$$
After the change of variables $s := K_Xt^{\al},\ t = K_X^{-1/\al}s^{1/\al},\ 
dt = K_X^{-1/\al}\al^{-1}s^{1/\al - 1}ds$
the integral $I(u)$ changes to
$$
C_X\IL_{K_Xu^{\al}}^{\infty}K_X^{-\be/\al}s^{\be/\al}K_X^{-1/\al}\al^{-1}s^{1/\al - 1}e^{-s}ds = C_XK_X^{-(\be+1)/\al}\al^{-1}\IL_{K_Xu^{\al}}^{\infty}s^{(\be + 1 - \al)/\al}e^{-s}ds.
$$
But for $\mu \in \Mr$ we have
$$
\IL_u^{\infty}s^{\mu}e^{-s}ds \backsim u^{\mu}e^{-u}.
$$
This is easily deduced from L'Hospital's rule:
$$
(u^{\mu}e^{-u})' = -u^{\mu}e^{-u} + \mu u^{\mu - 1}e^{-u} \backsim -u^{\mu}e^{-u} = \frac d{du}\IL_u^{\infty}s^{\mu}e^{-s}ds.
$$
Hence we easily obtain (4). The proof is complete. $\square$

\bigskip

Also, one can replace $\MP\{\zeta < -u\}$ by $\MP\{\zeta \le -u\}$, and similarly for $\eta$.
The asymptotic will remain the same.

\bigskip

{\bf Lemma 9.} Suppose $X$ is a random variable, $m := \ess\inf X$. Suppose $C_X, \al > 0$ are constants.

1. If $m  = -\infty$, then
$$
F_X(-u) = \MP\{X \le -u\} \backsim C_X(-u)^{-\al}\ \Leftrightarrow \ \MP\{X < -u\} \backsim C_X(-u)^{-\al}.
$$

2. If $m \in \Mr$, then, as $u \downarrow m$, we have:
$$
F_X(u) = \MP\{X \le u\} \backsim C_X(u - m)^{\al} \ \Leftrightarrow \ \MP\{X < u\} \backsim C_X(u - m)^{\al}.
$$

{\it Proof of Lemma 9.} Let us prove the first statement. For all $u$
$$
\MP\{X \le -u - 1\} \le \MP\{X < -u\} \le \MP\{X \le -u\} \le \MP\{X < -u + 1\}.
$$
If $\MP\{X \le -u\} \backsim C_X(-u)^{-\al}$, then 
$\MP\{X \le -u-1\} \backsim C_X(-u-1)^{-\al} \backsim C_X(-u)^{-\al}$, and $\MP\{X < -u\} \backsim C_X(-u)^{-\al}$. 
Similarly, if $\MP\{X < -u\} \backsim C_X(-u)^{-\al}$, then 
$\MP\{X < -u+1\} \backsim C_X(-u+1)^{-\al} \backsim C_X(-u)^{-\al}$, and
 $\MP\{X \le -u\} \backsim C_X(-u)^{-\al}$.

The proof of the second statement is similar; it is necessary to consider the inequalities
$$
\MP\{X \le u - u^2\} \le \MP\{X < u\} \le \MP\{X \le u\} \le \MP\{X < u + u^2\}.
$$
(We assume w.l.o.g. that $m = 0$.) $\square$

\bigskip

{\bf\Large 13. Auxillary lemmas.}

\bigskip

{\bf Lemma 10.} Suppose $a > 0$, $\mu > 0$, $f, S : [0, a] \to \Mr$ are continuous on $[0, a]$, 
$f(0) \ne 0$, $\min\limits_{[0, a]}S$ is attained only at the point $0$. Suppose $S \in C^1[0, \delta_0]$ 
for some $\delta_0 \in (0, a]$, $S'(0) > 0$. Then
$$
\CF(u) := \IL_0^ax^{\mu - 1}f(x)e^{-uS(x)}dx \backsim \CG(u) := \Ga(u)f(0)S'(0)^{-\mu}\la^{-\mu}e^{-uS(0)}.
$$

{\it Proof of Lemma 10.} W.l.o.g. suppose $f(0) \ne 0$. One can find $\delta_1 \in (0, \delta_0]$
such that $S'(x) > 0$ for all $x \in [0, \delta_1]$. 
Take an arbitrarily small $\eps > 0$ and find $\delta \in (0, \delta_1]$ such that
 $|f(x) - f(0)| < \eps f(0)$ for all $x \in (0, \delta]$. Then we have
$$
\CF(u) = \left(\IL_0^{\delta} + \IL_{\delta}^a\right)x^{\mu - 1}f(x)e^{-uS(x)}dx.
$$
The integral on $[\delta, a]$ is $O\left(e^{-Mu}\right)$, where $M := \min_{[\delta, a]}S > S(0)$ 
(Lemma 1.1, chapter 2, [3]), hence it is $o(\CG(u))$. And the integral on $[0, \delta]$ is estimated in this way:
$$
(1 - \eps)I(u) \le \IL_{\delta}^ax^{\mu - 1}f(x)e^{-uS(x)}dx \le (1 + \eps)I(u),\ I(u) := f(0)\IL_0^{\delta}x^{\mu - 1}e^{-uS(x)}dx.
$$
Change variables in $I(u)$: $S(x) - S(0) = S'(0)t,\ x = \hat{x}(t),\ \hat x \in C^1$. 
It maps $[0, \delta]$ into $[0, \delta']$ for some $\delta'$. Hence we obtain:
$$
I(u) = f(0)\IL_0^{\delta'}\hat x^{\mu-1}(t)e^{-u(S'(0)t + S(0))}\hat x'(t)dt.
$$
But $\hat x'(t) = S'(0)/S'(\hat x(t))$, $\hat x'(0) = 1$, hence the function 
$g(t) := \hat x^{\mu-1}(t)t^{-\mu+1}\hat{x}'(t), \ g(0) := 1$ is continuous on $[0, \delta']$. 
By Watson's lemma (see [2])
$$
I(u) = f(0)e^{-uS(0)}\IL_0^{\delta'}g(t)t^{\mu-1}e^{-uS'(0)t}dt \backsim f(0)(S'(0))^{-\mu}\Gamma(\mu)u^{-\mu}e^{-uS(0)} = \CG(u).
$$
Thus:
$$
1 - \eps \le \varliminf\limits_{u\to\infty}\frac{\CF(u)}{\CG(u)} \le \varlimsup\limits_{u\to\infty}\frac{\CF(u)}{\CG(u)} \le 1 + \eps.
$$
Since $\eps > 0$ is arbitrary, the proof is complete. $\square$

\bigskip

{\bf Lemma 11.} 1. Suppose $(X, \CM)$ is a measurable space with $\sigma$-finite measure $\mu$ on $\CM$. 
Let $f = f(u, x), g  = g(u, x): [a, \infty) \times X \to \Mr$ be measurable with respect to
$x \in \CM$ and Lebesgue integrable on a set $A \in \CM$ for every $u \in [a, \infty)$. 
If $f(x, u) \backsim g(x, u)$ uniformly for $x \in A$, then
$$
\IL_Xf(x, u)I_A(x)d\mu(x) \backsim \IL_Xg(x, u)I_A(x)d\mu(x).
$$
2. Suppose $(\Oa, \CF, \MP)$ is a probability space, $\xi_1 = (\xi_1(u), u \ge a), \xi_2 = (\xi_2(u), u \ge a)$
are random processes, and for the given event $A \in \CF$ $\xi_1(u)I_A$, $\xi_2(u)I_A$ 
are integrable for all $u \in [a, \infty)$. If $\xi_1(u) \backsim \xi_2(u)$ 
uniformly for $\omega \in A$, then
$$
\ME\xi_1(u)I_A \backsim \ME\xi_2(u)I_A.
$$
3. Suppose $I \subseteq \Mr$ is an interval, $f, g : [a, \infty) \times I \to \Mr$ are functions such that
for any $u \ge a$ they are Lebesgue-integrable on $I$ (and measurable by $t \in I$). 
If $t \in I$ $f(t, u) \backsim g(t, u)$ uniformly for $t \in I$, then
$$
\IL_If(t, u)dt \backsim \IL_Ig(u, t)dt.
$$
{\it Proof of Lemma 11.} Left to the reader. $\square$

\bigskip

{\bf\Large Acknowledgements.}

\bigskip

The author is deeply grateful to Prof. Vladimir Piterbarg, Lomonosov Moscow State University,
the Department of Mechanics and Mathematics, who supported me during this work and provided with valuable advice and comments. 

Also, the author is grateful to Prof. Krzysztof Burdzy and Prof. Soumik Pal (University of Washington, Seattle; the Department of Mathematics), who supported me when this article was being finished. I would like to thank them for useful discussion and important advice.

\bigskip

{\bf\Large References.}

\bigskip

[1] Jurg H$\ddot{u}$sler, V. I. Piterbarg. {\it Extremes of a certain class of Gaussian processes.}// Stochastic Processes and their Applications, N. 83, pp. 257 - 271. Elsevier, 1999.

[2] A. Erdelyi. {\it Asymptotic Expansions.} Dover, 1956.

[3] V. I. Piterbarg, {\it Asymptotical Methods in the Theory of Gaussian Processes and Fields.} AMS Translations of Mathematical Monographs 148, Providence, Rhode Island, 1996.

[4] Jurg H$\ddot{u}$sler, V. I. Piterbarg, Yueming Zhang. {\it Extremes of Gaussian Processes with Random Variance.} 2009.

[5]  Ioannis Karatzas, Steven E. Shreve. {\it Brownian Motion and Stochastic Calculus. Graduate Texts in Mathematics.} Springer, 2006.

[6] Norman L. Johnson, Samuel Kotz, N. Balakrishnan.  {\it Continuous univariate distributions. Wiley Series in Probability and Mathematical Statistics}, Vol. 1, 1994.

[7]  Ditlev Monrad, Holger Rootzen. {\it  Small values of Gaussian Processes and functional laws of the iterated logarithm.} Probability Theory Related Fields, 101 (1995), no. 2, 173 - 192.

[8] William Feller. {\it An Introduction to Probability Theory and Its Applications.} Vol. 1, 2nd edition, John Wiley \& Sons Inc., 1957. 

\end{document}